\documentclass[12pt]{article}
\usepackage{amsfonts}
\usepackage{mathrsfs, amssymb, amsmath, graphicx}
\usepackage{color}

\newenvironment{proof}[1][Proof]{\noindent\textbf{#1.} }{\ \rule{0.5em}{0.5em}}

\def\endproof{\nobreak\qquad\vrule height4pt width4ptdepth0pt}
\def\<{\langle}
\def\>{\rangle}
\def\AAA{\mathcal{A}}
\def\B{\mathcal{B}}
\def\A{\mathscr{A}}

\def\D{\mathcal{D}}
\def\G{\mathcal{G}}
\def\H{\mathcal{H}}
\def\K{\mathcal{K}}
\def\L{\mathcal{L}}

\def\N{\mathcal{N}}

\def\R{\mathcal{R}}

\def\S{\mathcal{S}}

\def\V{\mathcal{V}}

\def\phi{\varphi}

\begin{document}
\date{}
\baselineskip18pt
\title{Commutativity and Null Spaces of Unbounded Operators on Hilbert Space\footnote{R. Kadison
and S. Levin designed research, and wrote major part of the manuscript before Kadison's passing in 2018.
Levin and Z. Liu completed the remaining work and finished writing the article in 2022 and 2023.}}

\author{Richard Kadison\footnote{Department of Mathematics, University of Pennsylvania} \qquad Simon Levin\footnote{Department of Ecology and Evolutionary Biology, Princeton University}  \qquad
Zhe Liu\footnote{Department of Mathematics and Statistics, Georgetown University}
}
\maketitle
\begin{abstract}
Based on the success of a well-known method
for solving higher order linear differential equations, a study of
two of the most important mathematical features of that method, {\it viz.}
the null spaces and commutativity of the product of unbounded linear
operators on a Hilbert space is carried out.  A principle is proved
describing solutions for the product of such operators in terms of the
solutions for each of the factors when the null spaces of those factors
satisfy a certain geometric relation to one another.  Another geometric
principle equating commutativity of a closed densely defined operator and a
projection to stability of the range of the projection under the closed
operator is proved.
\end{abstract}


\section{Introduction}

One of the simplest and
most powerful steps underlying the solution of higher order linear
homogeneous differential equations with constant coefficients is the
factorization of the linear operator, allowing the representation of
solutions in terms of sums of solutions of the corresponding reduced
equations. For example, the solution of the equation
\vskip3pt
\centerline{$u''+ (\alpha+\beta)u'+\alpha\beta u=0,$}
\vskip3pt
\noindent where $\alpha$ and $\beta$ are distinct constants, is found by rewriting the
equation as
\vskip3pt
\centerline{$\big(\frac{d}{dx}+\alpha I\big)\big(\frac{d}{dx}+\beta I\big)u=0,$}
\vskip3pt
\noindent where $I$ is the identity operator; solving the equations
\vskip3pt
\centerline{$\big(\frac{d}{dx}+\alpha I\big)u=0  \ \ {\rm and} \ \ \big(\frac{d}{dx}+\beta I\big)u=0,$}
\vskip3pt
\noindent and then expressing the general solution as
\vskip3pt
\centerline{$u = ae^{-\alpha x}+be^{-\beta x},$}
\vskip3pt
\noindent where $a$ and $b$ are arbitrary constants, and the two terms in the general solution are respectively solutions of $(\frac{d}{dx}+\alpha I)u=0$ and $(\frac{d}{dx}+\beta I)u=0$.
This method breaks down if $\alpha = \beta$, since in this case the general solution collapses
to a one-parameter  family of solutions, insufficient to generate the full
set of solutions of the original equation. The simplest  example is, of course,
$u''=0$, whose solutions are $u=ax+b$.

Similar principles apply to partial differential equations of second order.
The most famous is the (2-dimensional) wave equation
\vskip3pt
\centerline{$\frac{\partial^2u}{\partial t^2}=c^2\frac{\partial^2u}{\partial z^2}$}
\vskip3pt
\noindent and we can change it into a simpler form
\vskip3pt
\centerline{$\frac{\partial^2u}{\partial xy}=0$}
\vskip3pt
\noindent by letting $x=z-ct$ and $y=z+ct$,
whose general solution may be written as (assuming continuous differentiability)
\vskip3pt
\centerline{$u=F(x)+G(y),$}
\vskip3pt
\noindent namely, as a sum of solutions of the factor equations
\vskip3pt
\centerline{$\frac{\partial u}{\partial x}=0  \ \ {\rm and} \ \ \frac{\partial u}{\partial y}=0.$}
\vskip3pt

As a final example, consider the simple system of linear equations
\vskip2pt
\centerline{$ABz=0,$}
\vskip3pt
\noindent where $A$ and $B$ are real commuting square matrices. Clearly, any solution of
\vskip2pt
\centerline{$Az=0  \ \ {\rm or} \ \  Bz=0$}
\vskip3pt
\noindent will also be a solution of $ABz=0$, and hence, so too will be any sum of such
solutions.  Conversely, the general solution of $ABz=0$ can be represented as a
direct sum of the solutions of $Az=0$ and $Bz=0$ provided their null spaces have
only the zero vector in common.  (This can be seen by ``counting.'' If $n_a$,
$n_b$, and $n_{ab}$ are the dimensions of the null spaces of $A$, $B$, and
$AB$, respectively, then the obvious fact that the null spaces of $A$ and $B$
are contained in the null space of $AB$ and the assumption of zero
intersection yield the inequality $n_a+n_b\le n_{ab}$.  At the same time, $B$
is a linear mapping of the null space of $AB$ into that of $A$, from $ABz=0$,
with kernel the null space of $B$.  Thus $n_{ab}-n_b\le n_a$, and
$n_{ab}=n_a+n_b$.)

These examples lead one to ask how general this method is.  In particular, if
$A$ and  $B$ are linear operators, under what conditions can the solutions
of $ABu=0$ be represented as sums of solutions of $Au=0$ and $ Bu=0$? The purpose of this
note is to address this question, with broad implications.  More
specifically, we seek rigorous mathematical statements and their proofs that
establish a relation between the null spaces of the factors in a product of
closed densely defined operators on a Hilbert space, the general framework
for differential operators and the $L_2$-solutions of their associated
differential equations, and the null space of the product.  In the next
section, we present some background information and preliminary material.  We
discuss the more involved concepts of ``commutativity'' and ``reducing'' in the
case of unbounded operators.  In Section 3, we begin with a general result
that provides a condition for the sum of the null spaces to have closure equal
to the null space of the closure of the product when that product is closable
and includes a condition for the sum of the null spaces to be closed.  We
apply this result to show that (unbounded) commuting self-adjoint operators
have null spaces whose sum is the null space of the closure of the product
(which is always closable).  As a corollary, the same is true for a
self-adjoint operator and a closed densely defined operator that commute.  We
then study the more general situation of a ``finite von Neumann algebra'' and
the operators ``affiliated'' with it.  We prove a ``rank-nullity'' theorem and
some ``relative-dimension'' inequalities for such operators and use these
results to show that the closure of the sum of those null spaces is the null
space of the closure of their product when those operators commute and their
null spaces are disjoint.  In Section 4, we study the consequences of
commutativity of a closed densely defined operator with a bounded
self-adjoint operator, especially with an orthogonal projection operator,
from a technical viewpoint.  In particular, we establish a geometric
principle equating commutativity of the closed operator and the projection
to the special stability, appropriate to these circumstances, which we
call ``proper stability,'' of the range of the projection under the closed
operator (Theorem 4.9).  The results in that section are closely
related to a larger program initiated by the
first and third-named authors ([1]).  That program undertakes to determine the structure of
the algebras of operators affiliated with finite von Neumann algebras.
Those affiliated algebras are called Murray-von Neumann algebras.
An important aspect of that program involves determining the nature of commuting
(self-adjoint) subfamilies of Murray-von Neumann algebras.  Later in Section 4, we
discuss, briefly, some of those results. (See [2], and [3] for more
commutativity results. Notably in [2], the
centers and the maximal abelian self-adjoint
subalgebras of Murray-von Neumann algebras are described.) In the last section, we tie our results to the
differential equations discussed in the introduction.

\section{Preliminaries}

In this section, we present some preliminary material that we
shall need in the following sections.  If $\H$ and $\K$ are Hilbert spaces (over the complex
numbers $\mathbb{C}$) and $\H\times\K$ is their Cartesian product, then $\H\times\K$
provided with the addition $(x_1,y_1)+(x_2,y_2)=(x_1+x_2,y_1+y_2)$ and
multiplication by scalars $a(x,y)=(ax,ay)$ is, again, a linear space (over
$\mathbb{C}$).  Provided, as well, with the inner product $\<(x_1,y_1),(x_2,y_2)\>=
\<x_1,x_2\>+\<y_1,y_2\>$, where $\<x_1,x_2\>$ and $\<y_1,y_2\>$ denote the
inner products of $x_1$ and $x_2$ in $\H$ and $y_1$ and $y_2$ in $\K$,
respectively, $\H\times\K$ with its linear space structure becomes a Hilbert
space, the {\it Hilbert direct sum\/} $\H\oplus\K$ of $\H$ and $\K$.

Suppose $T$ is a linear transformation from a linear submanifold $\D(T)$ (the
{\it domain\/} of $T$) of $\H$ with range in $\K$.  If $\D(T)$ is dense in $\H$,
we say that $T$ is {\it densely defined}.  The
{\it graph\/} $\G(T)$ of $T$ is the set $\{(x,Tx):\,x\in\D\}$ in
$\H\oplus\K$.  Since $T$ is linear, $\G(T)$ is a linear submanifold of
$\H\oplus\K$.  When $\G(T)$ is a closed subspace of $\H\oplus\K$, we say that
$T$ is {\it closed}. In addition, if $\H=\K$, we call $T$ a {\it closed operator\/} (on
$\H$).

Recall that an operator $B$ defined everywhere on $\H$ is continuous if and
only if it is {\it bounded\/} (that is, $\{\|Bx\|:\,x\in\H,\,\|x\|\le1\}$ is
a bounded set of numbers -- its least upper bound, denoted by $\|B\|$, is
called the {\it bound\/} or {\it norm\/} of $B$).  Each bounded operator $B$
has a closed graph (in $\H\oplus\H$).  To see that, note that if $(x,y)$ is
in the closure of $\G(B)$, there is a sequence $\{(x_n,Bx_n)\}$ in $\G(B)$
tending to $(x,y)$.  Thus $x_n\to x$ and $Bx_n\to y$.  Since $B$ is assumed
to be continuous (bounded), $Bx_n\to Bx$, and $Bx=y$.  Thus
$(x,y)=(x,Bx)\in\G(B)$, and $\G(B)$ is closed.  A famous and deep
(fundamental) result is a ``partial converse" to this: If $\D(T)$ is $\H$
($T$ is ``everywhere defined") and $\G(T)$ is closed, then
$T$ is bounded (continuous). This is the Closed Graph Theorem.

Although bounded operators are certainly the most convenient to handle,  many
of the operators we encounter (including most of the differential operators
and, of course, those discussed earlier) are unbounded.  That doesn't make
the case of bounded operators uninteresting or useless; the way we study the
unbounded operators we can handle is, largely, through the families of bounded
operators we can associate with them (through standard procedures -- for
example, ``spectral resolutions'').  The unbounded operators that are most
tractable are those with a closed graph.  That condition on the graph doesn't
replace continuity but it does substitute for it in some respects. An example of 
how the closed graph condition can substitute for
continuity (and one that is particularly important for what we
want to do) is the fact that the null space of a closed operator
$T$ is closed.  Suppose that $\{y_n\}$, in the null space of $T$, tends to $y$ in $\H$.  Then $\{(y_n,Ty_n)\}$ tends to $(y,0)$.  Since $\G(T)$ is closed,
$(y,0)\in\G(T)$.  In particular, $y\in\D(T)$ and $Ty=0$.  Thus
$y$ is in the null space of $T$, and that null space is closed.

If the graph of $T$ is not closed, it is natural to consider its closure
$\G(T)^-$ and ask if this closure is the graph of some operator
$\overline T$.  Now, it may happen that $\G(T)$ contains a sequence
$\{(x_n,Tx_n)\}$ such that $x_n\to 0$ and $Tx_n\to y$.  In that case,
$(0,y)\in\G(T)^-$.  If $y$ is not 0, there is no possibility for
$\G(T)^-$ to be the graph of a {\it linear\/} transformation
$\overline T$.  But if that doesn't occur, then, $\G(T)^-$ is the graph of a closed
linear transformation $\overline T$.  In that case, we say that $T$ is {\it
preclosed\/} (or {\it closable\/}) and that $\overline T$ is the {\it
closure\/} of $T$.  A preclosed (unbounded) operator often lends itself to
analysis by studying its closure.  Very often the ``natural'' unbounded
operators we encounter are, first, defined on some ``natural'' domain.  It is
only by enlarging that domain that we find the closure of the operator that
interests us, {\it if\/} it has a closure.

One reason it is so valuable to know that a densely defined operator $T$ is
closed is that, then, $T$ has a ``polar decomposition'': $T=VH$, where $V$ is
isometric on the range of $H$ and annihilates the orthogonal complement of
that range and $H$ ($=(T^*T)^\frac12$) is a
(positive) self-adjoint operator on $\H$ (generally unbounded when $T$ is).
The operator $H$ is subject to ``spectral analysis'' and $V$ is a fairly
simple bounded operator.

For a closed operator $T$, we say that a linear submanifold $\V$ of $\D(T)$
{\it reduces\/} $T$ when $\V$ and $\D(T)\ominus\V$ ($=(\H\ominus\V)\cap\D(T)$,
$\D(T)\ominus\V$ denotes the set of all vectors in $\D(T)$ that are orthogonal
to $\V$) generate a dense linear submanifold of $\D(T)$, $T(\V)\subseteq\V^-$, and
$T(\D(T)\ominus\V)\subseteq(\D(T)\ominus\V)^-$.

Commutativity for operators $A$ and $B$ in the purely algebraic sense is
simply the equality of $AB$ and $BA$.  In the case of operators on a Hilbert
space (or, more generally, a normed space), if $A$ and $B$ are bounded (and
everywhere defined), commutativity remains just the equality of $AB$ and $BA$.
When $A$ and $B$ are unbounded, even closed, this simple equality no longer
serves as the expression of an adequate concept of commutativity.  A hint of
what can cause difficulties can be seen by considering the case where $B$ is
the (everywhere defined) operator 0.  In this instance, $AB$ is, again, 0.
However, $BA$ is the 0 operator defined only on $\D(A)$, the domain of $A$.
Of course, we want to think of 0 as commuting with each closed operator.  We do
have that $A0$ is an {\it extension\/} of $0A$.  (We
say that $T$ is an extension of $S$ and write `$S\subseteq T$' when
$\D(S)\subseteq\D(T)$ and $Tx=Sx$ for each $x$ in $\D(S)$.)  For a closed
operator $A$ and a bounded operator $B$, the relation $BA\subseteq AB$ serves
as an adequate concept of ``commutativity'' (of $A$ and $B$).  In this case,
$AB$ is always closed, while $BA$, in general, will not be closed nor even
closable.  To see this, suppose $x_n\in\D(AB)$, $x_n\to x$, and $ABx_n\to y$.
Then note that $Bx_n\in\D(A)$ and $A(Bx_n)$
($=ABx_n$)$\to y$, and $Bx_n\to Bx$ since $B$ is continuous.  As $A$ is closed, $(Bx,y)\in\G(A)$ and $ABx=y$.  Thus
$(x,y)\in\G(AB)$, and $AB$ is closed.
An example showing that $BA$ need not be closed nor even closable can be
constructed by choosing $\{x\in\H\,:\,\sum_{n=1}^\infty n^4|\<x,e_n\>|^2<\infty\}$,
as the domain $\D$ in $\H$, where $\{e_1,e_2,\ldots\}$ is an orthonormal
basis for $\H$.  Let $B$ map $x$ to $\<x,z\>z$, for each $x$ in $\H$, where
$z$ is the vector $\sum_{n=1}^\infty n^{-1}e_n$.  (Not only is $B$ bounded,
but its range is one dimensional.)  Let $A$ map $x$ in $\D$ to
$\sum_{n=1}^\infty n^2\<x,e_n\>e_n$.  Since each finite linear combination
of the $e_n$ is in $\D$, $\D$ is dense in $\H$. To see that $A$ is closed,
suppose $x_m \in \D$, $x_m\to x$, and $Ax_m\to y$. For $e_k\in \{e_1,e_2,\ldots\}$,
$\<Ax_m, e_k\>=\<\sum_{n=1}^\infty n^2\<x_m,e_n\>e_n, e_k\>=k^2\<x_m, e_k\>\to k^2\<x, e_k\>$.
At the same time $\<Ax_m, e_k\>\to \<y, e_k\>$, so that $k^2\<x, e_k\>=\<y, e_k\>$
(for every $k$),
and $\sum_{n=1}^\infty n^4|\<x,e_n\>|^2=\sum_{n=1}^\infty |n^2\<x,e_n\>|^2=
\sum_{n=1}^\infty |\<y,e_n\>|^2=\|y\|^2< \infty$. It follows that $x\in \D$ and
$Ax=\sum_{n=1}^\infty n^2\<x,e_n\>e_n=\sum_{n=1}^\infty \<y,e_n\>e_n=y$.
Thus $A$ is closed. However, if $y_m=m^{-1}e_m$, then
$y_m\to 0$.  But $BAy_m=B\sum_{n=1}^\infty n^2\<y_m,e_n\>e_n=
B\sum_{n=1}^\infty n^2\<m^{-1}e_m,e_n\>e_n=mBe_m=m\<e_m,z\>z=
m\<e_m,\sum_{n=1}^\infty n^{-1}e_n\>z=z$ ($\ne 0$).  Thus $(0,z)\in\G(BA)^-$
and $BA$ is not closable.

If $A$ and $B$ are closed and densely defined, the question of their
commutativity is more delicate and complicated.  When $A$ and $B$ are
self-adjoint, the safest and most effective approach is to equate their
commutativity with that, of their spectral projections:
$E_\lambda F_{\lambda'}=F_{\lambda'}E_\lambda$, for all real $\lambda$ and
$\lambda'$, where $\{E_\lambda\}$ and $\{F_\lambda\}$ are the, respective,
spectral resolutions of $A$ and $B$.  In this case, $AB$ is closable and its
closure $\overline{AB}$ is self-adjoint.  Moreover, $A$ and $B$ are
``affiliated'' with the, respective, (abelian) ``von Neumann algebras''
generated by $\{E_\lambda\}$ and $\{F_\lambda\}$.  A {\it von Neumann algebra}
on $\H$ is an algebra of bounded operators on $\H$ that contains $I$,
contains $A^*$ when it contains $A$, and is strong-operator closed.  (The
strong-operator topology on $\B(\H)$ corresponds to the convergence,
$A_n\to A$ when $A_nx\to Ax$ for each $x$ in $\H$.)  The closed densely
defined operator $T$ is {\it affiliated} with the von Neumann algebra $\R$,
written as `$T\,\eta\,\R$',
when $U^*TU=T$ for each unitary operator $U$ that commutes with all operators
in $\R$.  (Note that the equality $U^*TU=T$ requires that $\D(U^*TU)=\D(T)$
and $U^*TUx=Tx$ for each $x$ in the common domain.)

The ``commutativity'' of $A$ and $B$, when they are assumed only to be closed
and densely defined, is a less useful and more obscure concept.  We might
identify this commutativity with affiliation of $A$ and $B$ with von Neumann
algebras $\R$ and $\S$, respectively, such that $RS=SR$ for each $R$ in
$\R$ and $S$ in $\S$.  But it is too strong.  It entails that $A$
commutes with $B^*$ as well as with $B$; this need not be the case even for
finite commuting matrices. When $A$ or $B$ is self-adjoint,
affiliation with commuting von Neumann algebras works well for commutativity
of $A$ and $B$.  If, say, $A$ is self-adjoint and $A$ and $B$ are affiliated
with von Neumann algebras such that the elements in one commute with all the
elements of the other, then $A$ commutes with the self-adjoint operator
$(B^*B)^{\frac12}$ and with the bounded operator $V$, the components of the
polar decomposition of $B(=V(B^*B)^{\frac12})$; that is to say, the projections in the spectral
resolutions of $A$ and $(B^*B)^{\frac12}$ commute with each other and
$VA\subseteq AV$.  The case of commutativity of bounded operators agrees with
what we have just defined.  If $AB=BA$ with $A$ self-adjoint, then
$B^*A=(AB)^*=(BA)^*=AB^*$.  It follows that each element in the von Neumann
algebra generated by $A$ commutes with each element in the von Neumann
algebra generated by $B$ and $B^*$.  Thus $A$ commutes with the components of
the polar decomposition of $B$ when it is self-adjoint and commutes with $B$.

This discussion does not provide us with a reasonable concept of
commutativity in the general case of closed densely defined operators $A$
and $B$.  With $VH$ and $WK$ the polar decompositions of $A$ and $B$, respectively, if we
were to say that $A$ and $B$ commute when $V$ and $H$ each commute with each of
$W$ and $K$, then $A$ might not commute with $A$, for that would require
that the components of the polar decomposition of $A$ commute with each
other.  That is equivalent to $A$ being ``normal'' ($A^*A=AA^*$), which is not
the case, in general, even for finite matrices.

For a special class of von Neumann algebras, those that are ``finite" in some
technical sense (to be explained), the affiliated operators behave very
nicely in relation to one another.  A linear functional $\tau$ on a
von Neumann algebra $\R$ that takes the value 1 at $I$ is said to be a
``trace" (on $\R$) when $\tau(AB)=\tau(BA)$ for each $A$ and $B$ in $\R$.
We say that $\R$ is a {\it finite} von Neumann algebra when it admits a
family of traces  such that $A=0$ if $\tau(A^*A)=0$ for each trace $\tau$ in
the family.  Of course, each bounded linear functional that is 1 at $I$ on an
abelian von Neumann algebra $\AAA$ is a trace on $\AAA$.  Thus from the
Hahn-Banach theorem, each abelian von Neumann algebra is finite in that
sense. An interesting class of finite von Neumann algebras is constructed
with the aid of discrete groups.  Let $G$ be such a group and let  $\H$ be
$l_2(G)$, the Hilbert space of absolutely square summable complex-valued
functions on $G$.  Let us also assume that $G$ is countably infinite so that
$\H$ is separable. If $x_g$ is the function that takes the value 1 at
the group element $g$ and $0$ at other elements of $G$, then $\{x_g:\,g\in G\}$ is an orthonormal basis for
$l_2(G)$. Left ``translations'' by $x_g$ with $g\in G$ ($=\{L_{x_g}: g\in G\}$) on the functions
in $l_2(G)$ give rise to unitary operators on $l_2(G)$ that generate an
algebra of operators whose strong-operator closure $\L(G)$ is a finite von
Neumann algebra (acting on $\H$) ([4]-II, Theorem 6.7.2, Proposition 6.7.4).
A simple computation shows that the linear functional
$A\to\<Ax_g,x_g\>$ (where $\<Ax_g,x_g\>$ is the inner product of the vectors
$Ax_g$ and $x_g$ in $\H$) is a trace on $\L(G)$ for each $g$ in $G$.  If
$A\in\L(G)$ and $\<A^*Ax_g,x_g\>=0$ for all $g$ in $G$, then $\|Ax_g\|^2=0$
for all such $g$ and $Ax_g=0$ for all such $g$.  As $\{x_g\}$ is an
orthonormal basis for $l_2(G)$, $A=0$.  Thus $\L(G)$ is finite.

If $\R$ is a finite von Neumann algebra (for instance, if $\R$ is abelian)
and $A$ and $B$ are affiliated with $\R$, then $AB$ is preclosed and
$\overline{AB}$ (denoted by $A\,\hat\cdot\,B$) is affiliated with $\R$, in particular,
$AB$ is densely defined ([4]-IV, Exercise 8.7.60).  Commutativity of $A$ and $B$, in this case, is the
equality of $A\,\hat\cdot\,B$ and $B\,\hat\cdot\,A$.
Since the algebra of $n\times n$ complex matrices is a finite von Neumann algebra
in the sense just defined, the concept of commutativity we have described
(by means of affiliation) coincides with the usual notion of commutativity.

We shall need a criterion for the sum $\mathfrak{X}+\mathfrak{Y}$
($=\{x+y:\,x\in\mathfrak{X},\,y\in\mathfrak{Y}\}$) of two closed disjoint subspaces $\mathfrak{X}$ and
$\mathfrak{Y}$ of a Banach space $\mathcal{Z}$ to be closed.  For this, we define the {\it
angle\/} between $\mathfrak{X}$ and $\mathfrak{Y}$ to be
$\inf\{\|x-y\|:\,x\in\mathfrak{X},\ y\in\mathfrak{X},\,\|x\|=\|y\|=1\}$.  It
is not difficult to prove that $\mathfrak{X}+\mathfrak{Y}$ is closed in $\mathcal{Z}$ if and only if
that angle is positive ([4]-I, pp. 63-64]).  Of course,
in the case of finite-dimensional spaces, all linear manifolds are closed and
this criterion plays no role there.

\section{Products}

In this section, we study the null space of
products of linear operators using the techniques discussed in the last
section.

\vskip6pt
{\bf Proposition 3.1.}\hskip8pt {\it Let $A$ and $B$ be closed linear
operators on a Hilbert space $\H$ with null spaces $\N_1$ and $\N_2$,
respectively.  Suppose $\V_1$ is a dense linear submanifold of $\N_1$ that
reduces $B$.  If $AB$ is closable, then
the null space $\N$ of $\overline{AB}$ is $[\N_1+\N_2]$, the closed
linear span of $\{u_1+u_2:\,u_1\in\N_1,\,u_2\in\N_2\}$.  If $\N_1$ and
$\N_2$ form a positive angle with one another, then $\N=\N_1+\N_2$.}
\vskip3pt
{\bf Proof.}\hskip4pt If $u_1\in\V_1$, then $Bu_1\in\N_1$, by assumption,
and $ABu_1=0$.  If $u_2\in\N_2$, then $Bu_2=0$, and $ABu_2=0$.  Thus $\N$
contains $\V_1+\N_2$.  As the null space of a closed operator is closed and
$\V_1$ is dense in $\N_1$, $\N$ contains $[\V_1+\N_2]$
($=[\N_1+\N_2]$).

If $u\in\N\ominus[\N_1+\N_2]$, then $u\in\D(AB)\subseteq\D(B)$, whence
$Bu\in\D(A)$ and $u$ is orthogonal to $\V_1$.  Thus $u\in\D(B)\ominus\V_1$.
By assumption, $Bu\in(\D(B)\ominus\V_1)^-$.  Hence $Bu\in\H\ominus\N_1$.
Assume that $u\ne0$.  Then, since $u\in\H\ominus\N_2$, $u$ is not in
$\N_2$ (otherwise, $\<u,u\>=\|u\|^2=0$, and $u=0$, contrary to assumption).
Thus $Bu\ne0$.  As $Bu\in\H\ominus\N_1$, $Bu$ is not in $\N_1$.  Since
$Bu\in\D(A)$, $ABu\ne0$, contrary to the choice of $u$ in $\N$.  Thus the
assumption that $u$ in $\N\ominus[\N_1+\N_2]$ is not 0 is untenable.
It follows that $\N=[\N_1+\N_2]$.
When $\N_1$ and $\N_2$ form a positive angle,
$\N_1+\N_2=[\N_1+\N_2]=\N$.\endproof
\vskip6pt

As noted in the preceding section, commuting unbounded self-adjoint
operators are particularly suited to joint study.  We apply the preceding
proposition to that case in the theorem that follows.
\vskip6pt
{\bf Theorem 3.2.}\hskip8pt{\it If $A$ and $B$ are commuting self-adjoint
operators acting on a Hilbert space $\H$ with null spaces $\N_A$ and $\N_B$, respectively, then
$\N_A+\N_B=\N_{A\,\hat\cdot\, B}$, where $\N_{A\,\hat\cdot\,B}$ is the null space
of $A\,\hat\cdot\,B$.}
\vskip3pt
{\bf Proof.}\hskip4pt Let $\AAA$ be the abelian von Neumann algebra generated
by $A$ and $B$ (note that $A$ and $B$ are not necessarily bounded here; see
Lemma 5.6.7 and Remark 5.6.11 of [4]-I for details).  Let $E$, $F$, and $G$, be the projections with respective
ranges $\N_A$, $\N_B$, and $\N_{A\hat\cdot B}$.  If $N$ is a projection
commuting with $A$, then $NA\subseteq AN$.  If $x\in\N_A$, then
$x\in\D(NA)\subseteq\D(AN)$ and $ANx=NAx=0$.  Hence $Nx\in\N_A$.  It follows
that, $NE=ENE=(ENE)^*=(NE)^*=EN$.  Thus $E$ is in the von Neumann algebra
generated by $A$, and hence, in $\AAA$.  Applying this conclusion to $B$ and
$A\,\hat\cdot\, B$, we have that $F$ and $G$ are also in $\AAA$.  Thus $E$, $F$, and
$G$, commute with one another and with $A$, $B$, and $A\,\hat\cdot\, B$.  In
particular, $EB\subseteq BE$.  If $x\in\D(B)$, then $Ex\in\D(B)$ and
$E(Bx)=BEx$.  Since $E$ and $B$ commute and are affiliated with the finite
(abelian) von Neumann algebra $\AAA$, $BE$ has a dense domain $\D$ (see
[4]-IV, Exercise 8.7.60]).  Let $\V=E(\D)$. Since $E$ is continuous, $\V$ is
dense in $E(\H)$.  We show that $\V$ reduces $B$, from which, with the aid
of the preceding proposition, $[\N_A+\N_B]=\N_{A\,\hat\cdot\, B}$.  To see this,
note that $\V\subseteq\D(B)$ since $\D$ is the domain of $BE$.  If $x\in\V$,
then $EBx=BEx=Bx$, whence $Bx\in E(\H)=\V^-$.  Similarly $I-E\in\AAA$ and $I-E$
commutes with $B$.  Again, $B(I-E)$ has a dense domain $\D'$ and $(I-E)(\D')$
($=\V'$) is dense in $(I-E)(\H)$.  Thus $\V'\subseteq\D(B)\ominus\V$.  Of
course, $\V+\V'$ is dense in $\H$ as is $\V+(\D(B)\ominus\V)$.  If
$y\in\D(B)\ominus\V$, then $y=(I-E)y$ since $y$ is orthogonal to $E(\H)$
($=\V^-$).  Thus $By=B(I-E)y=(I-E)By\in(I-E)(\H)=\V'^-=(\D(B)\ominus\V)^-$,
and $\V$ reduces $B$.

To complete the proof, we note that, since $E$ and $F$ commute, $\N_A+\N_B$
is closed, whence $\N_A+\N_B=\N_{A\,\hat\cdot\, B}$.  To see this, observe that
$E-EF$ and $F$ are orthogonal projections with sum a projection having range
$\N_A+\N_B$, which is, accordingly, closed.\endproof
\vskip6pt

Referring, again, to our discussion of commutativity, we had an acceptable
concept in the case where one of the operators, say $A$, is self-adjoint, and
$B$ is closed and densely defined.  We identified commutativity
of $A$ and $B$ with that of $A$ and $V$ and $A$ and $H$, where $VH$ is the polar decomposition of $B$.  In this case, $\N_B$ is the null
space of $H$.  From the preceding theorem, applied to $A$ and $H$, we have:
\vskip6pt

{\bf Corollary 3.3.}\hskip8pt{\it If $A$ is a self-adjoint operator and $B$ is a
closed densely defined operator acting on a Hilbert space $\H$ with null
spaces $\N_A$ and $\N_B$, respectively, then when $A$ and $B$ commute,
$\N_A+\N_B=\N_{A\,\hat\cdot\, B}$, where $\N_{A\,\hat\cdot\, B}$ is the null
space of $A\,\hat\cdot\, B$.}
\vskip6pt
Commutativity for arbitrary operators is not enough to guarantee that the
space generated by their null spaces is the null space of their product (even
of its closure).  This can be illustrated at the level of $2\times2$
matrices.  If $A$ is $({0\atop0}{1\atop0})$, then $A^2$ is
$({0\atop0}{0\atop0})$ with null space the full (two-dimensional) space, while
the null space of $A$ is not the full space.  Note, however, that the null
spaces of the ``two'' commuting operators ($A$ and $A$) generate a subspace of
the null space of their product ($A^2$).  This was remarked on in the introduction,
where we also observed that if the null spaces of the commuting operators are
disjoint, then the closed linear manifold they generate is the null space of
the product.  We shall prove these results for operators affiliated with a
finite von Neumann algebra.  The center-valued dimension function $\Delta$
on the set of all projections in the finite von Neumann algebra
([4]-II, Section 8.4]) is needed for that discussion.

We begin with an extension of the ``rank-nullity'' theorem to the case of
operators affiliated with a finite von Neumann algebra. For the theorems that follow, we shall use the fact, noted in the preceding
section, that $TB$ is closed when $T$ is closed and $B$ is bounded.
We use $N(T)$ and $R(T)$ to denote the projections whose ranges are,
respectively, the null space of $T$ and the closure of the range of $T$.

\vskip6pt
{\bf Theorem 3.4.}\hskip8pt{\it If $\R$ is a finite von Neumann algebra with
center-valued dimension function $\Delta$ and $T\,\eta\,\R$, then
\vskip3pt
\centerline{$\Delta(R(T))+\Delta(N(T))=I.$}
\vskip3pt
\noindent If $E$ and $F$ are projections in $\R$ such that $FTE=TE$, then
$\Delta(E)-\Delta(N(T)\wedge E)\le\Delta(F).
$
In particular, if $N(T)\wedge E=0$, then
$\Delta(E)\le\Delta(F)$.}
\vskip3pt
{\bf Proof.}\hskip4pt From Exercise 6.9.52 of [4]-IV, $R(T)\sim R(T^*)$. From Exercise 2.8.45 of [4]-III, $N(T)=I-R(T^*)$.  Thus
$\Delta(R(T))+\Delta(N(T))=\Delta(R(T^*))+\Delta(N(T))=\Delta(I)=I$.

Since $N(TE)=N(T)\wedge E+I-E$,
$\Delta(N(TE))=\Delta(N(T)\wedge E)+\Delta(I-E)$.
Hence
\begin{align}
I&=\Delta(R(TE))+\Delta(N(TE))\nonumber\\
&=\Delta(R(TE))+\Delta(N(T)\wedge E)+\Delta(I-E),\ \text{and}\nonumber\\
\Delta(E)&=\Delta(I)-\Delta(I-E)=I-\Delta(I-E)\nonumber\\
&=\Delta(R(TE))+\Delta(N(T)\wedge E).\nonumber
\end{align}
Now, $R(TE)\le F$, since $FTE=TE$, by assumption.  Thus
$\Delta(E)-\Delta(N(T)\wedge E)=\Delta(R(TE))\le\Delta(F)$.
If $N(T)\wedge E=0$, then
$\Delta(E)\le\Delta(F)$.\endproof

\vskip6pt
{\bf Theorem 3.5.}\hskip8pt{\it If $\R$ is a finite von Neumann algebra with
center-valued dimension function $\Delta$, and $A$ and $B$ are operators
affiliated with $\R$ such that $A\,\hat\cdot\,B=B\,\hat\cdot\,A$, then
$[\N_A+\N_B]\subseteq\N_{A\,\hat\cdot\, B}$.
If $\N_A\cap\N_B=\{0\}$, then $[\N_A+\N_B]=\N_{A\,\hat\cdot\, B}$.}
\vskip3pt

{\bf Proof.}\hskip4pt If $x\in\N_A$, then $Ax=0\in\D(B)$, whence
$x\in\D(BA)\subseteq\D(B\,\hat\cdot\, A)$, and
$A\,\hat\cdot\,Bx=B\,\hat\cdot\,Ax=BAx=0.$
Hence $x\in\N_{A\,\hat\cdot\, B}$, and $\N_A\subseteq\N_{A\,\hat\cdot\, B}$.
Symmetrically, $\N_B\subseteq\N_{A\,\hat\cdot\, B}$.  As $\N_{A\,\hat\cdot\, B}$ is
closed (in the Hilbert space $\H$ on which $\R$ acts), and
$[\N_A+\N_B]\subseteq\N_{A\,\hat\cdot\, B}$.

Assume, now, that $\N_A\cap\N_B=\{0\}$ so that $\Delta(N(A)\wedge N(B))=0$.
From the formula in Exercise 8.7.31 of [4]-IV,
\begin{align}
&\Delta(N(A))+\Delta(N(B))\nonumber\\
=&\Delta(N(A)\vee N(B))+\Delta(N(A)\wedge N(B))\nonumber\\
=&\Delta(N(A)\vee N(B)).\nonumber
\end{align}
As we noted, $[\N_A+\N_B]\subseteq\N_{A\,\hat\cdot\, B}$.  Thus
\vskip3pt
\centerline{$\Delta(N(A))+\Delta(N(B))=\Delta(N(A)\vee N(B))\le\Delta(N(A\,\hat\cdot\, B)),$}
\vskip3pt
\noindent and $N(A)\vee N(B)\le N(A\,\,\hat\cdot\,B)$. We shall establish the reverse inequality.

Let $WK$ be the polar decomposition of $B$ and $F_n$ be the spectral
projection for $K$ corresponding to the interval $[-n,n]$ for each positive
integer $n$.  Let $G_n$ be the projection $F_n\wedge(N(A\,\hat\cdot\, B)-N(B))$.
Then $BF_n$ is everywhere defined and bounded, whence
\vskip3pt
\centerline{$(A\,\hat\cdot\, B)\,\hat\cdot\, F_n= A\,\hat\cdot\,(B\,\hat\cdot\, F_n)=A\,\hat\cdot\,(BF_n)=ABF_n,$}
\vskip3pt
\noindent since $ABF_n$ is closed, as noted earlier.  If $x\in G_n(\H)$, then
$x\in F_n(\H)\subseteq\D(B)$ and
$x\in\N(A\,\hat\cdot\, B)\subseteq\D(A\,\hat\cdot\, B)$. Thus
$0=(A\,\hat\cdot\, B)x=(A\,\hat\cdot\, B)F_nx=ABF_nx=ABx.$
Hence $Bx\in\N_A$.  When $x\ne0$, $Bx\ne0$, since $x$ is orthogonal to
$\N_B$.  Thus $B$ is a one-to-one mapping of $G_n(\H)$ into $\N_A$.
Algebraically, $N(B)\wedge G_n=0$ and $N(A)BG_n=BG_n$.  From the preceding
theorem, $\Delta(G_n)\le\Delta(N(A))$.  Since $F_n\uparrow I$ and $\R$ is
finite, from Exercise 8.7.34(ii) of [4]-IV,
$G_n\uparrow(N(A\,\hat\cdot\, B)-N(B))$.
It follows from Exercise 8.7.33(i) of [4]-IV that
$\Delta(G_n)\uparrow\Delta(N(A\,\hat\cdot\,B)-N(B))
=\Delta(N(A\,\hat\cdot\, B))-\Delta(N(B))$.
Hence $\Delta(N(A\,\hat\cdot\,B))-\Delta(N(B))\le\Delta(N(A))$.  When
$\Delta(N(A)\wedge N(B))=0$,
$\Delta(N(A\,\hat\cdot\,B))\le\Delta(N(A))+\Delta(N(B))=\Delta(N(A)\vee N(B))
+\Delta(N(A)\wedge N(B))=\Delta(N(A)\vee N(B))$.
Since $N(A)\vee N(B)\le N(A\,\hat\cdot\,B)$, we conclude that
$[\N_A+\N_B]=\N_{A\,\hat\cdot\,B}$ when we assume that
$\N_A\cap\N_B=\{0\}$.  If the angle between $\N_A$ and $\N_B$ is
positive, then $\N_A+\N_B=\N_{A\,\hat\cdot\,B}$.\endproof

\section{Commutativity and stability}

When $A$ is a closed densely
defined operator on a Hilbert space $\H$ and $B$ is a bounded operator on
$\H$, we say that $A$ and $B$ commute when $BA\subseteq AB$.  Under these
assumptions on $A$ and $B$, this turns out to be a useful definition of
commutativity, especially because the assumption of ``closedness'' is present.
Some of the consequences of this commutativity are proved in this section.

\vskip6pt
{\bf Proposition 4.1.}\hskip8pt{\it If $B$ and $C$ are, respectively, a
bounded self-adjoint and a closed densely defined operator, acting on a
Hilbert space, for which $BC\subseteq CB$, then $BC^*\subseteq C^*B$.}
\vskip3pt
{\bf Proof.}\hskip4pt Since $BC\subseteq CB$, $(CB)^*\subseteq(BC)^*$ from Remark 2.7.5 of [4]-I.  From Lemma 6.1.10 of [4]-II, $(BC)^*=C^*B^*=C^*B$.  From Exercise 2.8.44 of [4]-III, and the fact that $CB$ is densely defined (as
$BC\subseteq CB$ and $BC$ is densely defined with $\D(BC)=\D(C)$), we have
$BC^*=B^*C^*\subseteq(CB)^*$.  Thus
$BC^*\subseteq(CB)^*\subseteq(BC)^*=C^*B. \endproof$
\vskip6pt
{\bf Corollary 4.2.}\hskip8pt {\it If $E$ is an orthogonal projection operator and
$C$ is a closed densely defined operator acting on a Hilbert space, and
$EC\subseteq CE$, then $EC^*\subseteq C^*E$.}
\vskip6pt
The preceding results tell us that if a closed densely defined operator $C$
commutes with a bounded self-adjoint operator, then so does its adjoint.  If
$C$ is bounded as well as $B$, then this conclusion is a simple consequence
of the equality
$BC^*=B^*C^*=(CB)^*=(BC)^*=C^*B^*=C^*B$,
when $CB=BC$.  As we see, this becomes a somewhat more complicated situation
when an unbounded $C$ is involved.  At the same time, the preceding corollary
takes us very close to the geometric interpretation of an operator $C$
commuting with a projection $E$.  In case $C$ is bounded that commutativity
is equivalent to the range of $E$ being ``stable'' (``invariant'') under $C$
and $C^*$.  If $EC=CE$ and $x\in E(\H)$, then $Cx=CEx=ECx\in E(\H)$.  This
same conclusion is valid for $C^*x$, for as just noted, $EC^*=C^*E$ when
$EC=CE$.  For the converse, note that $E(\H)$ is stable under $C$ if and only
if $CE=ECE$.  To see this, observe that if $C$ maps $E(\H)$ into itself, then
$CEx\in E(\H)$ for each $x$ in $\H$, whence $CEx=ECEx$ and $CE=ECE$.
On the other hand, if $CE=ECE$ and $x\in E(\H)$, then $Cx=CEx=ECEx\in E(\H)$,
and $E(\H)$ is stable under $C$.  It follows, now, that if $E(\H)$ is stable
under both $C$ and $C^*$, then $CE=ECE$ and $C^*E=EC^*E$, whence $EC=ECE=CE$,
taking adjoints of both sides of the equality $C^*E=EC^*E$.

Of course, one wonders, at this point, how the familiar principle whose proof we
have just given, the stability of $E(\H)$ under the bounded operator $C$ and
its adjoint is equivalent to commutativity of $E$ and $C$, fares when we drop
the assumption that $C$ is bounded and everywhere defined.  At the very
beginning, we must adjust what we mean by stability under $C$, since $C$ acts
only on $\D(C)$.  We can't expect $C$ to map all of $E(\H)$ into $E(\H)$.  In
this case, stability under $C$ must mean that $C$ maps $\D(C)\cap E(\H)$ into
$E(\H)$.  Commutativity of $E$ and $C$, now means that $EC\subseteq CE$, of
course.  From this commutativity, $CE$ must be densely defined since $\D(C)(=\D(EC)$) is contained in $\D(CE)$.  Thus $E(\D(C))\subseteq\D(C)$,
and, therefore, $E(\D(C))\subseteq\D(C)\cap E(\H)$.  At the same time,
$E(\D(C))$ is dense in $E(\H)$, as $\D(C)$ is dense in $\H$ and $E$ is
continuous on $\H$.  The conclusion that $C$ maps $\D(C)\cap E(\H)$ into
$E(\H)$ when $E$ and $C$ commute is significant since, in this case,
$E(\D(C))\subseteq\D(C)\cap E(\H)$, and $E(\D(C))$ is dense in $E(\H)$.  In
Theorem 4.9, we demonstrate that the ``operator-projection commutativity
principle'' holds not only when the operator is bounded, but also when it is
closed and densely defined. For the purposes
of this theorem, we define ``proper stability'' of $E(\H)$ under $C$ with the
observations we have just made.
\vskip6pt
{\bf Definition 4.3.}\hskip8pt The range $E(\H)$ of an orthogonal projection
operator $E$ on a Hilbert space $\H$ is said to be {\it properly stable\/}
under a closed densely defined operator $C$ on $\H$ when $C$ maps
$\D(C)\cap E(\H)$ into $E(\H)$ and $E(\D(C))\subseteq\D(C)$.
\vskip6pt
{\bf Remark 4.4.}\hskip8pt If $E(\D(C))\subseteq\D(C)$, then $E(\D(C))\subseteq\D(C)\cap E(\H)$.  However,
with $y$ in $\D(C)\cap E(\H)$, $y=Ey\in E(\D(C)$, whence
$\D(C)\cap E(\H)\subseteq E(\D(C))$, and $E(\D(C))=\D(C)\cap E(\H)$, under
the assumption that $E(\D(C))\subseteq\D(C)$.
\vskip6pt
{\bf Lemma 4.5.}\hskip8pt {\it If $C$ is a closed densely defined operator and $E$ is an orthogonal projection operator on a Hilbert space $\H$, then $C$ maps $\D(C)\cap E(\H)$ into $E(\H)$
if and only if $CE=ECE$.}
\vskip3pt
{\bf Proof.}\hskip4pt  Note that $\D(CE)=\D(ECE)$. It will suffice to show
that one of $CE$ or $ECE$ is an extension of the other to show that they are
equal.  Suppose that $C$ maps $\D(C)\cap E(\H)$ into $E(\H)$ and that
$x\in\D(CE)$.  Then $Ex\in\D(C)\cap E(\H)$.  By assumption,
$CEx\in E(\H)$, whence $ECEx=CEx$.  Since this holds for each $x$ in $\D(CE)$,
$CE\subseteq ECE$, and $CE=ECE$.

Suppose, now, that $CE=ECE$.  If $x\in\D(C)\cap E(\H)$, then $Ex=x\in\D(C)$ so
that $x\in\D(CE)$ and $x\in\D(ECE)$.  Moreover, $Cx=CEx=ECEx\in E(\H)$.  Thus
$C$ maps $\D(C)\cap E(\H)$ into $E(\H)$ in this case. \endproof
\vskip6pt
{\bf Corollary 4.6.}\hskip8pt{\it With the notation of the preceding lemma, if $C$ maps $\D(C)\cap E(\H)$ into $E(\H)$, then $ECE$ is closed.}
\vskip3pt
{\bf Proof.}\hskip4pt From the preceding lemma, $CE=ECE$. Then $CE$ is closed since $C$ is closed and $E$ is bounded. \endproof
\vskip6pt
{\bf Lemma 4.7.}\hskip8pt{\it If $C$ is a closed densely defined operator and
$E$ is an orthogonal projection operator on a Hilbert space $\H$ such that
$CE$ and $C^*E$ are densely defined, and $C$ and $C^*$ map, respectively,
$\D(C)\cap E(\H)$ and $\D(C^*)\cap E(\H)$ into $E(\H)$, then
$CE\subseteq\overline{EC}$ and $C^*E\subseteq\overline{EC^*}$.  In addition,
$EC$ and $EC^*$ are densely defined and preclosed.}
\vskip3pt
{\bf Proof.}\hskip4pt From Lemma 4.5, $CE=ECE$ and $C^*E=EC^*E$.  Since $CE$
and $C^*E$ are closed, so are $ECE$ and $EC^*E$. From Theorem 2.7.8 of [4]-I,
$ECE=(ECE)^{**}$ and $EC^*E=(EC^*E)^{**}$.  From Exercise 2.8.44 of [4]-III, $EC^*\subseteq(CE)^*$
since $CE$ is densely defined, by assumption.  Hence $EC^*$ is
preclosed (with dense domain $\D(C^*)$), and similarly for $EC$.  Thus, from
Lemma 6.1.10 of [4]-II and fact 5.6(13) of [4]-I,
$(EC)^*=C^*E=EC^*E\subseteq(CE)^*E=(ECE)^*E$.
It follows from Remark 2.7.5 and Theorem 2.7.8 of [4]-I that $((ECE)^*E)^*\subseteq(EC)^{**}=\overline{EC}$. Now, again, from Exercise 2.8.44 of [4]-III,
$E(ECE)^{**}\subseteq((ECE)^*E)^*\subseteq\overline{EC}$.
Since $ECE=E(ECE)=E(ECE)^{**}\subseteq\overline{EC}$,
$CE=ECE\subseteq\overline{EC}$.
Arguing symmetrically, $C^*E\subseteq\overline{EC^*}$. \endproof
\vskip6pt
{\bf Lemma 4.8.}\hskip8pt{\it If $C$ is a closed densely defined operator and
$E$ is an orthogonal projection operator on a Hilbert space $\H$ such that
$C$ and $E$ commute, then $CE=\overline{EC}$ and $C^*E=\overline{EC^*}$.}
\vskip3pt
{\bf Proof.}\hskip4pt Since $C$ and $E$ commute, $EC\subseteq CE$.  As discussed earlier, $C$ maps $\D(C)\cap E(\H)$ into $E(\H)$
and $CE$ is closed and densely defined.  From Corollary 4.2, $C^*$ and $E$
commute; hence $C^*$ maps $\D(C^*)\cap E(\H)$ into $E(\H)$  and $C^*E$ is
closed and densely defined.  Thus the preceding lemma applies and
$CE\subseteq\overline{EC}$. As $CE$
is closed and $EC\subseteq CE$, it follows that $\overline{EC}\subseteq
CE\subseteq\overline{EC}$.  Hence $\overline{EC}=CE$.  Applying this
conclusion, with $C^*$ in place of $C$,  $\overline{EC^*}=C^*E$. \endproof
\vskip6pt

{\bf Theorem 4.9.} (Operator-Projection Commutativity Principle)\hskip8pt
{\it A closed densely defined operator $C$ and an orthogonal projection
operator $E$ on a Hilbert space $\H$ commute, that is, $EC\subseteq CE$, if and only if, the range of
$E$ is properly stable under $C$ and $C^*$, in which case, $C^*$ and $E$ also
commute.}
\vskip3pt
{\bf Proof.}\hskip4pt Note that the conditions for $E(\H)$ to be properly stable under
$C$ follow from the commutativity of $E$ and $C$.  This same commutativity
implies that $E$ and $C^*$ commute (Corollary 4.2). Then $E(\H)$ is properly stable under $C^*$.

We suppose, now, that $E(\H)$ is properly stable under $C$ and $C^*$.  Since
$E(\D(C))\subseteq\D(C)$, $\D(C)\subseteq\D(CE)$, whence $CE$ and, similarly,
$C^*E$ are densely defined and closed. From Lemma 4.7, $EC$ and $EC^*$ are
densely defined and preclosed.  Moreover, from Lemma 4.5, $CE=ECE$ and $C^*E=EC^*E$.  Our
present goal is to show that $EC\subseteq CE$.

If $x\in\D(EC)$, then $x\in\D(C)$ and $Ex\in\D(C)$ (since $E(\D(C))\subseteq\D(C)$).
Thus $x\in\D(CE)$ and $\D(EC)\subseteq\D(CE)$. To show that $EC\subseteq CE$,
we shall show that $ECx = CEx$ for each $x\in \D(EC)$. Now, $Ex\in\D(EC)(=\D(C))$ and
\vskip3pt
\centerline{$EC(Ex)=ECE(Ex)=CE(Ex)=CEx.$}
\vskip3pt
\noindent Note, too, that $(I-E)x=x-Ex\in\D(C)$.
It follows that $(I-E)x\in\D(EC)$.  Of course, $(I-E)x\in\D(CE)$
since $E(I-E)x=0\in\D(C)$.  From Lemma 4.7, $CE\subseteq\overline{EC}$.  It
follows that $(I-E)x\in\D(\overline{EC})$ and $\overline{EC}(I-E)x=CE(I-E)x=0$.
However, $(I-E)x\in\D(EC)$ as noted, and $EC\subseteq\overline{EC}$.
Thus
\vskip3pt
\centerline{$EC(I-E)x=\overline{EC}(I-E)x=CE(I-E)x=0.$}
\vskip3pt
\noindent Now, for each $x$ in $\D(EC)$, $x=Ex+(I-E)x$, and we have
\vskip3pt
\centerline{$ECx=ECEx+EC(I-E)x=CEx+0=CEx.$}
\vskip3pt
\noindent Hence $EC\subseteq CE$, and $C$ commutes with $E$. \endproof
\vskip6pt

As discussed earlier, the essence of commutativity, wherever it is seen
and to whichever mathematical system it applies, is embodied in the simple
algebraic expression $AB=BA$.  We saw some of the
difficulties that can occur when studying unbounded operators on a Hilbert
space $\H$, where very careful attention must be paid to the domains of the
operators that appear at each stage of an argument.  In addition, the
fact that the operators appearing at each stage are closed or closable must be monitored.  In the end, for an effective use of what
we can press from the properties of unbounded operators, we must ``descend''
from ``equality'' (``='') to ``extension'' (``$\subseteq$''), for
the most part.  Even then, what we saw in the preceding results of this
section was a study of commutativity of a bounded operator $B$ and a closed
densely defined operator $C$, expressed as $BC\subseteq CB$.  Is there
anything that can be said about a question as basic as the commutativity of
two unbounded self-adjoint operators $A$ and $B$?  Is there any vestige of
the original expression $AB=BA$, perhaps, in slightly modified form, from
which we can draw information, when $A$ and $B$ are self-adjoint?  There is
such a possibility, but before describing that, let us note some of the
limitations to the question we are answering when it is posed as broadly as
we have formulated it.  The operator $AB$, with $A$ and $B$ self-adjoint may
have only 0 in its domain -- even when $B$ is bounded.  As an example of
this, let $A$ be self-adjoint with domain $\D(A)$ a dense subset different
from $\H$ (so, $A$ is not bounded), let $x_0$ be a unit vector not in
$\D(A)$, and let $B$ be the orthogonal projection operator whose range is
spanned by the vector $x_0$.  Then $AB$ has $\{0\}$ as its domain, while
$\D(BA)=\D(A)$, dense in $\H$.  Of course,
$AB\subseteq BA$, indeed, $AB\subseteq T$ for each operator $T$ on
$\H$.  There is
clearly, not much ``nutrition'' to be had from considering operators such as
$AB$.  We must require of $AB$ and $BA$ that they be densely defined.
Considering commuting, bounded, self-adjoint operators $A$ and $B$, we see
that $AB=BA=B^*A^*=(AB)^*$.  Thus $AB$ is self-adjoint, in this case, if and only if $A$ and $B$ commute.  We can't expect $AB$ to be
closed, in general.  It may not be preclosed,
even if it is densely defined.  However, self-adjoint operators are closed;
it seems prudent to require that $AB$ and $BA$ are densely defined and
preclosed as well.  If commutativity of $A$ and $B$ is to retain some hint of
its most basic meaning, ``$AB=BA$,'' then whatever this equality is to mean in
that expression, it should be the case that $ABx=BAx$ whenever
$x\in\D(AB)\cap\D(BA)$($=\D$).  As we shall see, when the spectral
resolutions of $A$ and $B$ commute with one another (our strongest and most
``reassuring'' form of commutativity for $A$ and $B$), $\D$ is a core for both
$\overline{AB}$ and $\overline{BA}$ on which they agree.

We recall that a {\it core\/} for a closed densely defined
operator $T$ is a linear submanifold $\D_0$ of the domain $\D(T)$ of $T$,
dense in $\H$, such that
$\G(T|\D_0)^-=\G(T)$.
Submanifolds of $\D(T)$ with this graph property of $\D_0$ are often very
effective for working with $T$, especially when $\D_0$ consists of elements
in $\D(T)$ on which the action of $T$ is particularly simple and well
understood.  Of course, each linear submanifold of $\D(T)$ containing $\D_0$
is a core for $T$.  If $S$ is another closed densely defined operator and
$\D_0$ is a core for each of $S$ and $T$ on which $S$ and $T$ agree, then
$S|\D_0=T|\D_0$, whence $\G(S)=G(S|\D_0)^-=G(T|\D_0)^-=G(T)$, and $S=T$.  
If $S$ and $T$ are densely defined but just preclosed and $\D_0$ is a core for both $\overline{S}$ and
$\overline{T}$ on which they agree, then we can conclude that
$\overline{S}=\overline{T}$, as just proved.
If we know, in addition, that $\D_0\subseteq\D(S)\cap\D(T)(=\D)$, then $\D$
is a core for $\overline{S}$ and $\overline{T}$, since $\D_0$ is assumed to
be a core for each of them.  Thus $\overline{S|\D}=\overline{T|\D}$, since
$\overline{S}=\overline{T}$, as we proved.  Now, $S|\D\subseteq \overline{S|\D}$, and
$T|\D\subseteq \overline{T|\D}$.
As $S$ and $T$ are not assumed to be closed, we cannot proceed from the fact
that $S$ and $T$ agree on the intersection $\D$ of their domains to the
conclusion that $S$ and $T$ are equal; $\D(S)$ and $\D(T)$ can be different.

Addressing the question we raised earlier, in the light of the information
gathered in the preceding discussion, we shall prove that
$\overline{HK}=\overline{KH}$ for self-adjoint $H$ and $K$, under the
appropriate closure and density assumptions, as discussed, if and only if the
spectral resolutions of $H$ and $K$ commute (Theorems 4.17 and 4.18).  As we
noted, when $S$ and $T$ agree on a core for both of $\overline{S}$ and
$\overline{T}$, they agree on $\D(S)\cap\D(T)$, which is, in this case,
a core for each of $\overline{S}$ and $\overline{T}$.  It is precisely this
form of ``agreement'' for $\overline{HK}$ and $\overline{KH}$ that we assume
to argue that the spectral resolutions for $H$ and $K$ commute. We present a series of related commutativity results in the rest of this section. 
\vskip6pt

{\bf Lemma 4.10.}\hskip8pt{\it If $C$ is a closed densely defined operator
 and $B$ is a bounded everywhere defined operator on a Hilbert space $\H$,
then $C+B$ is closed with domain $\D(C)$.}
\vskip3pt
{\bf Proof.}\hskip4pt If $(x_n,y_n)\in\G(C+B)$, $x_n\to x$, and $y_n\to y$,
in $\H$, then $Cx_n=y_n-Bx_n\to y-Bx$, since $B$ is continuous.  As $C$ is
closed, $x\in\D(C)$ and $Cx=y-Bx$.  Thus $(x,y)\in\G(C+B)$, and $C+B$ is
closed. \endproof
\vskip6pt

The special case of this lemma that we shall use is the fact that $H+aI$ (and
$K+bI$) is closed when $H$ (and $K$) are self-adjoint and $a$ (and $b$) is
some complex number.
\vskip6pt
{\bf Lemma 4.11}\hskip8pt{\it If $C$ is a preclosed densely defined operator and $B$ is
a bounded everywhere defined operator on a Hilbert space $\H$, then $C+B$ is preclosed
and $\overline{C+B}=\overline{C}+B$.}
\vskip3pt

{\bf Proof.}\hskip4pt Suppose that $(0,y)\in\G(C+B)^-$ and $\{x_n\}$ is a
sequence in $\D(C+B)$($=\D(C)$) such that $x_n\to 0$ and $(C+B)x_n\to y$.  In
this case, since $B$ is continuous, $Cx_n\to y$.  As $C$ is preclosed, $y=0$
and $C+B$ is preclosed.

To see that $\overline{C+B}=\overline{C}+B$, observe that $\overline{C}+B$ is
closed, from Lemma 4.10, and $C+B\subseteq\overline{C}+B$.  Thus
$\overline{C+B}\subseteq\overline{C}+B$.  On the other hand, if
$x\in\D(\overline{C}+B)=\D(\overline{C})$, then there is a sequence $\{x_n\}$
in $\D(C)$ such that $x_n\to x$ and $Cx_n\to\overline{C}x$.  Since $B$ is
continuous, $Bx_n\to Bx$, whence $Cx_n+Bx_n=(C+B)x_n\to\overline{C}x+Bx$.
Note that $x_n\in\D(C)=\D(C+B)$.  Thus
$(x,\overline{C}x+Bx)\in\G(C+B)^-=\G(\overline{C+B})$.  In particular,
$x\in\D(\overline{C+B})$, and $\overline{C}+B\subseteq\overline{C+B}$.  It
follows that  $\overline{C}+B=\overline{C+B}$. \endproof 

\vskip6pt

If $C$ is preclosed and densely defined, then $C+B$ has
a densely defined adjoint ([4]-I, Theorem 2.7.8(ii)).
\vskip6pt
{\bf Lemma 4.12.}\hskip8pt{\it If $C$ is a preclosed and densely defined operator and $B$ is
a bounded everywhere defined operator on a Hilbert space $\H$, then $(C+B)^*=C^*+B^*$.}
\vskip3pt

{\bf Proof.}\hskip4pt Note that $y\in\D(C^*+B^*)$ ($=\D(C^*)$) if and only if
the mapping $x\to\<Cx,y\>$ is a bounded linear functional on $\D(C)$, and
this holds if and only if there is a $y'$ ($=C^*y$) such that
$\<Cx,y\>=\<x,y'\>$ for all $x$ in $\D(C)$.  Since $x\to\<Bx,y\>=\<x,B^*y\>$
is also a bounded linear functional on $\D(C)$, $x\to\<(C+B)x,y\>$ is a
bounded linear functional on $\D(C)$($=\D(C+B)$) if and only if $x\to\<Cx,y\>$
is; that is, $y\in\D((C+B)^*)$ if and only if $y\in\D(C^*)$($=\D(C^*+B^*$)).
When $y\in\D(C^*)$,
$$
\<x,(C+B)^*y\>=\<(C+B)x,y\>
=\<x,C^*y\>+\<x,B^*y\>=\<x,(C^*+B^*)y\>
$$
for all $x$ in $\D(C)$.  Since $\D(C)$ is (assumed to be) dense in $\H$,
$(C+B)^*y=(C^*+B^*)y$.  Hence $(C+B)^*=C^*+B^*$. \endproof
\vskip6pt
Next lemma can be formulated in purely algebraic, vector-space form.  The
spaces may be finite or infinite-dimensional.
\vskip6pt

{\bf Lemma 4.13.}\hskip8pt{\it If $H$ and $K$ are linear transformations
defined on domains $\D(H)$ and $\D(K)$, linear submanifolds of a vector
space $\V$, $A$ and $B$ are linear transformations with domains $\V$,
$A(\D)\subseteq\D(K)$ and $B(\D)\subseteq\D(H)$, where $\D=\D(HK)\cap\D(KH)$,
then
$\D(HK)\cap\D(KH)=\D((H+A)(K+B))\cap\D((K+B)(H+A))$.}
\vskip3pt
{\bf Proof.}\hskip4pt Note that $\D(HK)=(K^{-1}(\D(H))\cap\D(K)$, since
$x\in\D(HK)$ if and only if $x\in\D(K)$ and $Kx\in \D(H)$, which is the case,
if and only if $x\in\D(K)$ and $x\in K^{-1}(\D(H))$.  Thus
$$\D(HK)\cap\D(KH)=\D(H)\cap\D(K)\cap(H^{-1}(\D(K))\cap(K^{-1}(\D(H))),
$$
and since $\D(H+A)=\D(H)$ and $\D(K+B)=\D(K)$,
\begin{align}
&\D((H+A)(K+B))\cap\D((K+B)(H+A))\nonumber\\
=&\D(H+A)\cap\D(K+B)\cap((H+A)^{-1}(\D(K+B)))\cap((K+B)^{-1}(\D(H+A)))\nonumber\\
=&\D(H)\cap\D(K)\cap((H+A)^{-1}(\D(K)))\cap((K+B)^{-1}(\D(H)))\nonumber\\
=&\D(H)\cap\D(K)\cap(H^{-1}(\D(K)))\cap(K^{-1}(\D(H)))\nonumber\\
=&\D(HK)\cap\D(KH).\nonumber
\end{align}
The third inequality uses the fact that, if $x\in\D$, then
$x\in(H+A)^{-1}(\D(K))$, that is $(H+A)x\in\D(K)$ if and only if $Hx\in\D(K)$
(since $Ax\in A(\D)\subseteq\D(K)$), and the symmetric fact with $H$ and
$K$ interchanged and $A$ and $B$ interchanged. \endproof
\vskip6pt

{\bf Remark 4.14.}\hskip8pt Replacing $A$ by $aI$ and $B$ by $bI$, where $I$
is the identity transformation on $\V$, we have
$\D((H+aI)(K+bI))\cap\D((K+bI)(H+aI))=\D(HK)\cap\D(KH)$.
It is this last equality that we shall use in the proof of Theorem 4.18, with
$a$ and $b$ having the values of $\pm i$, to argue from what we prove for
unbounded (normal) operators to the commutativity of associated bounded
(normal) operators, and thence, to the abelian von Neumann algebras they
generate (via Theorems 2.7.8, 5.6.15, 5.6.18, Lemmas 5.6.7, 5.6.10,
Proposition 2.7.10, and Remark 2.7.11 of [4]-I).  From the vantage point of
the powerful results of [5] as broadened and strengthened in [6] and [7], the
proof of Theorem 4.17 is completed using the Murray-von Neumann algebra of
operators affiliated with this abelian von Neumann algebra.
\vskip6pt

{\bf Lemma 4.15.}\hskip8pt{\it If $H$ and $K$ are self-adjoint operators and
$A$ and $B$ are bounded everywhere defined operators on a Hilbert space $\H$, $\D=\D(HK)\cap\D(KH)$ is dense in $\H$ with $A(\D)\subseteq\D(K)$, $B(\D)\subseteq\D(H)$,
$A^*(\D)\subseteq\D(K)$, and $B^*(\D)\subseteq\D(H)$, then $(H+A)(K+B)$ and
$(K+B)(H+A)$ are preclosed with $\D$ in the domain of each of them.}
\vskip3pt

{\bf Proof.}\hskip4pt From Lemma 4.13, $\D=\D((H+A)(K+B))\cap\D((K+B)(H+A))$,
whence $\D\subseteq\D((H+A)(K+B))$, $\D\subseteq\D((K+B)(H+A))$, and each of
$H$, $K$, $(H+A)(K+B)$, and $(K+B)(H+A)$ is densely defined.

We show, next, that $(H+A)(K+B)$ and $(K+B)(H+A)$ are preclosed.  From Lemma 4.10, $H+A$ and $K+B$
are closed and densely defined.  From Lemma 4.12, $(H+A)^*=H+A^*$ and
$(K+B)^*=K+B^*$. Thus, with $x$ in $\D((H+A)(K+B))$ and $y$ in $\D(=\D(HK)\cap\D(KH))$,
$y\in\D(H)=\D(H+A^*)=\D((H+A)^*)$ and
$$
(H+A)^*y=Hy+A^*y\in\D(K)=\D(K+B^*)=\D((K+B)^*),
$$
so that the mapping
$$
x\to\<(H+A)(K+B)x,y\>=\<(K+B)x,(H+A^*)y\>=\<x,(K+B^*)(H+A^*)y\>
$$
is a bounded linear functional on $\D((H+A)(K+B))$ and
$y\in\D(((H+A)(K+B))^*)$. Hence $\D\subseteq\D(((H+A)(K+B))^*)$ and $((H+A)(K+B))^*$ is densely defined.
It follows that $(H+A)(K+B)$ (and, similarly, $(K+B)(H+A)$) is preclosed ([4]-I, Theorem 2.7.8(ii)).\endproof
\vskip6pt
{\bf Lemma 4.16.}\hskip8pt{\it If $C$ is a preclosed densely defined operator on a Hilbert space $\H$ and the range of $\overline{C}$ is dense in
$\H$, $\D$ is a core for $\overline{C}$, and $\D\subseteq\D(C)$, then $C(\D)$
is dense in $\H$.}
\vskip3pt
{\bf Proof.}\hskip4pt Since $\D$ is a core for $\overline{C}$,
$\G(C|\D)^-=\G(\overline{C})=\G(C)^-$.  Hence, with $y$ in the range of
$\overline{C}$, there is a vector $x$ in $\D(\overline{C})$ such that
$(x,y)\in\G(\overline{C})$ and some $(x',y')$ in $\G(C|\D)$ with $x'$ near $x$
and $y'$ near $y$.  Now, $y'\in C(\D)$. Thus $C(\D)$ is dense in the range of
$\overline{C}$, which is dense in $\H$.  It follows that $C(\D)$ is dense in
$\H$. \endproof
\vskip6pt

{\bf Theorem 4.17}\hskip8pt{\it Let $H$ and $K$ be self-adjoint operators
on a Hilbert space $\H$.  Let $\{E_\lambda\}_{\lambda\in\mathbb{R}}$ and
$\{F_\lambda\}_{\lambda\in\mathbb{R}}$ be the respective spectral resolutions of $H$
and $K$.  If $E_\lambda F_{\lambda'}=F_{\lambda'}E_\lambda$ for all real
$\lambda$ and $\lambda'$, then $\{E_\lambda\}_{\lambda\in\R}$,
$\{F_\lambda\}_{\lambda\in\R}$ generate an abelian von Neumann algebra $\AAA$
with which $H$ and $K$ are affiliated as are $H+iI$, $H-iI$, $K+iI$, $K-iI$,
with respective domains $\D(H)$, $\D(H)$, $\D(K)$, $\D(K)$.  Each of these
operators is an injective linear transformation of its domain onto $\H$ with
respective inverses $T_+$, $T_-$, $S_+$, and $S_-$, which are bounded, and
which generate $\AAA$.  All the intersections
$$
\D((H+aI)(K+bI))\cap\D((K+bI)(H+aI))
$$
are the same dense linear submanifold $\D$ of $\H$ for all choices of
complex scalars $a$ and $b$.  In particular, $\D=\D(HK)\cap\D(KH)$.  In
addition,
$(H+aI)(K+bI)x=(K+bI)(H+aI)x$,
for each $x$ in $\D$, and $\D$ is a core for $\overline{(H+aI)(K+bI)}$ and
$\overline{(K+bI)(H+aI)}$.  The products $HK$ and $KH$ are preclosed,
$$
\overline{HK}=H\,\hat\cdot\,K=K\,\hat\cdot\,H=\overline{KH},
$$
and $\overline{HK}$ is self-adjoint.}
\vskip3pt

{\bf Proof.}\hskip4pt Referring to Section 2, the ``safest and
most effective'' way to impose commutativity on unbounded self-adjoint
operators on $\H$ is to require that their spectral resolutions commute, the
condition $E_\lambda F_{\lambda'}=F_{\lambda'}E_\lambda$, which we are
presently assuming.  As we noted there, $\{E_\lambda,F_\lambda\}_{\lambda\in\R}$
generates an abelian von Neumann algebra $\AAA$ with which $H$ and $K$ are
affiliated.  Moreover, $\overline{HK}=\overline{KH}$.  That is,
$H\,\hat\cdot\,K=K\,\hat\cdot\,H$ in the algebra $\A_f(\AAA)$, the ``Murray-von
Neumann algebra,'' of operators affiliated with $\AAA$ ([5], [6], [7], Theorem 5.6.15 of [4]-I), and as argued in the proof of Theorem 5.6.15 of
[4]-I, $\D_0=\bigcup_{n=1}^\infty E^{(n)}F^{(n)}(\H)$, where $E^{(n)}=E_n-E_{-n}$
and $F^{(n)}=F_n-F_{-n}$,
is a core for $\overline{(H+aI)(K+bI)}$ and $\overline{(K+bI)(H+aI)}$, where $H+aI$ and $K+bI$
are closed, by Lemma 4.10. It follows that $\D$ (containing $\D_0$) is a core for $\overline{(H+aI)(K+bI)}$ and $\overline{(K+bI)(H+aI)}$. The statement about equality of domains
follows from Lemmas 4.13 and 4.15.  In addition,
$(H\,\hat\cdot\,K)^*=K^*\,\hat\cdot\,H^*=K\,\hat\cdot\,H$. Thus
$$
(\overline{HK})^*=(H\,\hat\cdot\,K)^*=K\,\hat\cdot\,H=\overline{KH}=\overline{HK},
$$
and $\overline{HK}$ is self-adjoint.  The assertions involving $T_+$, $T_-$,
$S_+$, and $S_-$, follow from Proposition 2.7.10 and Lemma 5.6.7 of [4]-I.
\endproof
\vskip6pt

The arguments alluded to in the preceding proof require a good deal of
spectral theoretic and functional analysis for their complete and careful
presentation.  In Sections 5.2 and 5.6 of [4]-I, they receive such a
presentation.  As a consequence, we allow the sketch following the statement
of the preceding theorem to suffice, together with detailed references.

In the theorem that follows, we make use of the notation of the preceding
theorem while abandoning the assumption of commuting spectral resolutions,
the strong form of commutativity for the self-adjoint operators $H$ and $K$,
as initial hypothesis and replacing it by the ``simple-minded'' (``classical'')
assumption of agreement of $HK$ and $KH$ on a large enough domain.  The goal
of that theorem is to recapture the ``commuting-resolutions'' condition.
\vskip6pt
Note, from Lemma 4.15, that $(H+aI)(K+bI)$ and $(K+bI)(H+aI)$, with $a$ and
$b$ in $\{-i,i\}$, are preclosed, for which reason, the assumption that each
of their closures has $\D$ as a core makes sense.
\vskip6pt

{\bf Theorem 4.18.}\hskip8pt{\it With the notation of the preceding theorem,
if $HK$ and $KH$ agree on $\D$ and $\D$ is a core for each of
$\overline{(H+aI)(K+bI)}$, where $a$ and $b$ are in $\{i,-i\}$, then
$E_\lambda F_{\lambda'}=F_{\lambda'}E_\lambda$ for all real $\lambda$ and
$\lambda'$, $H$ and $K$ are affiliated with $\AAA$, which is abelian, and
$\overline{HK}$ is self-adjoint.}
\vskip3pt

{\bf Proof.}\hskip4pt  Our plan is to show that $\{T_+,T_-,S_+,S_-\}$ is a
commuting family.  Each of $\{T_+,T_-\}$, $\{S_+,S_-\}$ generates an abelian
von Neumann algebra with which $H$ and $K$ are, respectively, affiliated.
Those abelian von Neumann algebras commute with each other and generate an
abelian von Neumann algebra $\AAA$ that contains both of them and, accordingly,
contains the spectral resolutions $\{E_\lambda\}_{\lambda\in\R}$ and
$\{F_\lambda\}_{\lambda\in\R}$.  Thus these resolutions commute with one
another. That algebra $\AAA$ is generated by $\{E_\lambda\}_{\lambda\in\R}$ and
$\{F_\lambda\}_{\lambda\in\R}$, as well; $H$ and $K$ are affiliated with $\AAA$,
and the remaining assertions follow from this.

We prove, first, that $((H+aI)(K+bI))(\D)$ is a dense linear submanifold
of $\H$ for any $a$ and $b$ in $\{-i,i\}$.  Since $H+aI$ and $K+bI$ are
injective mappings of their domains onto $\H$, the same is true for their
products.  (Recall, for this, that each element of the domain of $H+aI$ is
the image of some element of the domain of $K+bI$, which element is, by
definition, an element of the domain of the product $(H+aI)(K+bI)$.) Thus
Lemma 4.16 applies, and all the images  $((H+aI)(K+bI))(\D)$ (and
$((K+aI)(H+bI))(\D)$) are dense in $\H$.  Letting $T_a$ be $T_+$ if $a$ is $i$
and $T_-$ if $a$ is $-i$, $S_b$ be $S_+$ if $b$ is $i$ and $S_-$ if $b$ is
$-i$, we show that $T_aS_b=S_bT_a$ for the four possible products $T_aS_b$.
Note, for this, that if $x\in\D$ and $y=(H+aI)(K+bI)x$, then
$$
S_bT_ay=S_bT_a(H+aI)(K+bI)x=S_b(K+bI)x=x.
$$
Since
$y=(H+aI)(K+bI)x=HKx+bHx+aKx+abx
=KHx+aKx+bHx+bax=(K+bI)(H+aI)x$,
$$
T_aS_by=T_aS_b(K+bI)(H+aI)x=x.
$$
Thus $S_bT_ay=T_aS_by$, for each $y$ in $(H+aI)(K+bI)(\D)$, a dense subset of
$\H$ (as we have just proved).  As $S_bT_a$ and $T_aS_b$ are continuous
everywhere defined mappings of $\H$ into itself, $S_bT_a=T_aS_b$. \endproof
\vskip6pt

{\bf Remark 4.19.}\hskip4pt It is worth noting that the products $(H+aI)(K+bI)$
with $a$ and $b$ in $\{-i,i\}$ are actually closed.  We are in a position to complete the proof of Theorem
4.18 using only the information that those products are preclosed, which we
note by specializing, without further effort, from Lemma 4.15.  When $a$ and
$b$ take on the values $i$ and $-i$, as we saw, $(H+aI)(K+bI)$ is an
injective mapping of its domain onto $\H$, and the mapping $S_bT_a$ is its
set-theoretic inverse.  Thus $\phi(\G(S_bT_a))=\G((H+aI)(K+bI))$, where
$\phi$ is the homeomorphism of $\H\times\H$, in its product topology, onto
itself that transforms $(x,y)$ in $\H\times\H$ to $(y,x)$.  It follows that
one graph is closed if the other is closed.  However $S_bT_a$, as an
everywhere-defined, bounded operator, has a closed graph.  Thus $(H+aI)(K+bI)$
is closed.

\section{Back to basics}

Our ``Introduction'' consists of a few
``rough-and-ready'' illustrations of the factoring method for solving
differential equations in a large class of linear differential equations.
There is no mention in that discussion of what a solution is, to what family
of functions the equation is supposed to apply, and in what sense the formal
algebraic ``factorization'' of the equation represents the same problem as the
original equation.  The purpose of the sections following the introduction
is to make a set of choices among the many that could be made, consistent
with some of the uses and practices, and then, to provide (what we hope will
be seen as) an absolutely rigorous mathematical framework for the
factorization technique, noting at the same time, the extent to which that
technique is valid.  The first choice is to limit our study to
$L_2$-functions on spaces where differentiation has a clear and well-known
meaning, the real line for ordinary differential equations and $\mathbb{R}^n$ for
partial differential equations.  With these choices, we have the vast theory
of Hilbert spaces to apply to the function families, and the well-developed
machinery of linear operators and algebras of such operators on those Hilbert
spaces to apply to the differential equations.  It is, without doubt, the
case that the clever algebraic-combinatorial tricks one devises in
conjunction with the construction and study of ``special functions'' are at the
heart of the study and solving differential equations.  Whether or not we
recognize the facts that the exponential function gives us eigenvectors for
the differentiation operator or that the sine and cosine functions generate a
stable (invariant), two-dimensional subspace for that operator, in those
broader (linear-space) terms, we had better bring them to some high level of
awareness in our consciousness before we venture forward with the study of
differential equations.  We must also recognize that the most important force
behind the desire (the {\it need\/}) to study and solve differential
equations is their applications to astronomy, physics, engineering, chemistry,
biology, medicine, sociology, finance, and let us not forget, mathematics
itself -- just about everything that claims some relation to ``science.''
The formulation of the laws governing the structures to which these
applications are being made, within the appropriately chosen mathematical
model, are often formulated in terms of differential equations.  A prime
example of this is the Newton-Hamilton formulation of Newton's laws for the
(time) evolution of (classical) dynamical systems ([8], [9]).  In this case, as in so
many of the applications, the central goal is to obtain {\it numbers\/} from
calculations with solutions, which are, then, checked against well-designed,
carefully executed, associated experiments in the laboratory or against
earlier accumulated data.

A paradigm for this process is Bohr's ([10], [11]) theoretical calculation of the wave
lengths of the bright lines in the visible range (the ``Balmer series'') of
the spectrum of hydrogen.  Using Rutherford's ``planetary'' model of the
hydrogen atom, Bohr notes that the possible orbits of the single
(``planetary'') electron have radii such that the energy emitted by the
electron in dropping from one ``shell'' to another with smaller radius, at a
given frequency $\nu$ is an integral multiple of $h\nu$, as Planck was forced
to assume when deriving his formula for a full (black body) radiator (where
$h$ is Planck's experimentally determined ``quantum of action'' $\approx
6.625\times10^{-27}$ erg sec).  The calculations resulted in wave lengths
that agreed with accurately measured, experimental, spectroscopic data to a
remarkable degree of precision.

One may well wonder whether it is worth the effort and thought needed to put
the clever methods one devises, for  solving the equations and problems, that
yield the numbers, then verified with such accuracy in the laboratory, in a
rigorous mathematical framework.  Does anyone ({\it should anyone}) care, when
the methods yield the numbers and the laboratories are available for
verifying them?  We do not ask this question as an exercise in aesthetics:
Is an interest in the careful mathematical formulation, when the numbers and
laboratories are so effectively at hand, healthy or morbid?  That aspect
would seem to be one of ``individual nature'' and ``individual interest.''  It
is intended as a question of whether there is anything useful, of practical
importance, that emerges from such a ``rigorization'' exercise.  If we
formulate that question in terms of Bohr's paradigm, and note what has
occurred since that extraordinary success, we arrive at an answer that seems
to have broad applicability to our question.  Translated in terms of that
paradigm, our question becomes that of whether it is worth attempting to find
a firm mathematical basis for the inspired ``{\it ad hoc}'' quantum
assumptions such as those made by Planck, Einstein and Bohr (though, Planck's
assumption, the first of these, marking the beginning of Quantum Mechanics,
was engendered more by ``despairing submission'' than ``inspiration'').  Having
the history of how Quantum Mechanics developed following Bohr's discovery,
we can easily answer our question affirmatively.  The path that Dirac ([12])
takes is nicely illustrative of what seeking more than the numbers and
experimental verification can yield.  Dirac notes that having a classical
dynamical system analogous to the ``small'' system being studied can be an
extremely important feature of the problem, and performing classical
analytic-dynamical calculations, bolstered by an inspired (``{\it ad hoc}'')
choice of quantum assumption on that system, can lead to convincing agreement
with experimental data in the quantum problem.  He then recognizes that,
since the process of forming the ``Poisson bracket'' of two functions
plays so crucial a role in Hamilton's formulation of classical dynamics that
it must have a counterpart, with fundamental importance, in Quantum
Mechanics.  Reasoning informally and proceeding with cautious steps, Dirac
concludes that, whatever algebraic structure we fasten to replace the algebra
of functions in classical dynamics, ``Lie bracketing,'' $[A,B]=AB-BA$, should
replace Poisson bracketing.  From there, and by comparison with what Poisson
bracketing produces when applied to the ``canonical coordinates'' and their
``conjugate momenta'' of a classical dynamical system, Dirac extracts the
Heisenberg relation, the fundamental relation of Quantum Mechanics.  Its form
($[Q,P]=$)$QP-PQ=i\hbar I$ announces that whatever mathematical model we
use for Quantum Mechanics, its algebraic structure must be non-commutative
since $\hbar$ ($=\frac h{2\pi}$) is non-zero.  The full panoply of {\it ad
hoc\/} quantum assumptions is cleverly packed into this non-commutativity
assumption by the insertion of $\hbar$ as the basic amount of
non-commutativity needed.

At the same time, the Heisenberg relation forbids us from using the algebras
of finite matrices of various orders as our full model for the mathematical
description of Quantum Mechanics.  The trace functional, no matter how it is
normalized, is not 0 on $i\hbar I$ and is 0 on $QP-PQ$.  The algebra of all
matrices of a given order might otherwise be the model of choice for a system
with a finite number of particles.  Further, purely mathematical investigation
shows that even the (von Neumann) algebra of all bounded operators on an
infinite-dimensional Hilbert space will not support the Heisenberg relation.

It is one of those fortunate scientific events, that occur often enough to
free us from refereeing to them as ``miracles,'' that the mathematics suitable
for representing the Heisenberg relations and providing the model that Dirac
[12] needed had just been developed, notably by von Neumann ([13], [14]) and Marshall
Stone ([15], [16]) while they both were in their mid to late twenties.  Our
reference is to the theory of Hilbert spaces, the algebras of bounded linear
operators acting on them, and families of unbounded linear operators on them.
So, the unbounded operators on Hilbert spaces that we have been studying in
connection with the factorization method for solving a class of differential
equations appears, here, in the basic mathematical formulation of Quantum
Mechanics.  It is not surprising, then, that one of the operators appearing
in the classical representation of the Heisenberg relation on Hilbert space
is differentiation itself, suitably adorned and ``prepared'' for this role.
(See [1] for an extensive account of this classical representation
as well as a detailed discussion of limitations to Hilbert space
representations of the relation, some old and some new limitations.)

With reference to the question we posed earlier: Is there a point to
introducing mathematical rigor when informal methods produce results in close
agreement with otherwise-obtained data -- the Einstein-Bohr paradigms,
coupled with the Dirac-Heisenberg-Schr\"odinger mathematical formulations
give us the simple answer, ``Yes.''  Dirac moves forward, introducing the
operator algebras we now call abelian ``von Neumann algebras'' ([12]).
These algebras are key to what we have been doing in the earlier sections and
will reappear in a crucial way, shortly, when we tie that earlier work to the
examples in the introduction.  From this process of formalizing the
approximate mathematical model, Dirac predicts the positron, Schr\"odinger ([17])
finds his wave equation and states of dynamical systems evolving in time, and
Heisenberg ([18]) describes dynamical systems in terms of the expectations of
observables evolving in time.  In this instance, it would seem that there is
some value in the process of mathematical formulation and rigorization.

Turning, now, to the work at hand and casting the process of solving
differential equations in a Hilbert-space, unbounded-operator framework, in
particular, when factorization is possible, we have examined the workings of
that in a mathematically rigorous manner.  It remains to tie the results of
that examination to equations such as those studied, loosely, in the
introduction.  Our first task is to say precisely what we mean when we refer
to differentiation as a linear operator -- on which space? -- with which
domain? -- with which action?  To attempt this for the purposes of solving
differential equations is too large a project with too many disparate goals
and requirements.  The needs of quantum physics, however, with which we
illustrated the point to a careful mathematical formulation, give us one
clear and eminently worthwhile direction in which to head, one for which the
preceding sections have prepared us.  We shall realize differentiation as a
linear operator on Hilbert space, viewed as an $L_2$-space.  Differentiation
is elegantly displayed, through the medium of M. H. Stone's result [St2] as
the ``generator'' of a (continuous) one-parameter, unitary group, $t\to U_t$,
where $(U_t(f))(s)=f(s+t)$ for $f \in L_2(\mathbb{R})$ and all real $s$ and $t$.  For the reader's
convenience, we state Stone's theorem.
\vskip6pt
{\bf Stone's Theorem.}\hskip4pt{\it If $H$ is a (possibly unbounded)
self-adjoint operator on the Hilbert space $\H$, then $t\to\exp itH$ is a
one-parameter unitary group on $\H$.  Conversely, if $t\to U_t$ is a
one-parameter unitary group on $\H$, there is a (possibly unbounded)
self-adjoint operator $H$ on $\H$ such that $U_t=\exp itH$ for each real $t$.
The domain of $H$ consists, precisely, of those vectors $x$ in $\H$ for which
$t^{-1}(U_tx-x)$ tends to a limit as $t$ tends to 0, in which case, this
limit is $iHx$.}
\vskip6pt
A close examination of which equivalence classes of functions in $L_2(\mathbb{R})$ lie
in the domain of $H$ when $t\to U_t$ is the one-parameter unitary group on
$L_2(\mathbb{R})$ arising from letting $\mathbb{R}$ act by translation on $\mathbb{R}$ reveals that
they are those classes each of which contains a (necessarily unique)
function $f$ that is absolutely continuous on $\mathbb{R}$ with almost everywhere
derivative $g$ that lies in $L_2(\mathbb{R})$.  In this case,
$\frac1t(e^{itH}f-f)\to g$ as $ t\to0$,
where the limit is convergence in $L_2(\mathbb{R})$ (convergence ``in the mean of
order 2''), and $iH(\{f\})=\{g\}$, where $\{f\}$ and $\{g\}$ are the
equivalence classes of $f$ and $g$, respectively.  The fact that $g$ is the
derivative, almost-everywhere, of $f$ is the very strong motivation for
identifying $iH$, from the one-parameter, translation, unitary group,
$t\to\exp itH$, as the {\it differentiation\/} operator on $L_2(\mathbb{R})$(=$\H$)
with the domain and action as described.

Now, we tie related results obtained in this paper to the differential equations shown in
the introduction in a Hilbert-space, unbounded-operator framework (as discussed earlier, we limit our study to
$L_2$-functions for now). First, consider the second order ordinary differential equation example:
\vskip3pt
\centerline{$u''+ (\alpha+\beta)u'+\alpha\beta u=0,$}
\vskip3pt
\noindent with $\alpha$ and $\beta$ distinct constants, factored as
\vskip3pt
\centerline{$\big(\frac{d}{dx}+\alpha I\big)\big(\frac{d}{dx}+\beta I\big)u=0.$}
\vskip3pt
\noindent As just discussed, from Stone's Theorem,
we treat $\frac{d}{dx}$
as $iH$, which is a skew-adjoint operator defined on a certain domain in $\H$. As we know, the
operator $H$ is self-adjoint and it is then affiliated with some abelian von Neumann algebra
$\R$ (which is a finite von Neumann algebra). Of course, the operators
\vskip3pt
\centerline{$A=\frac{d}{dx}+\alpha I \ \ {\text{and}} \ \ B=\frac{d}{dx}+\beta I$}
\vskip3pt
\noindent are affiliated with $\R$ as well, and they commute. From Theorem 3.5,
$[\N_A+\N_B]\subseteq\N_{A\,\hat\cdot\, B}$, and
$[\N_A+\N_B]=\N_{A\,\hat\cdot\, B}$ when $\N_A\cap\N_B=\{0\}$,
which is the present case (as $\alpha$ and $\beta$ are distinct).
Now, consider the wave equation example:
\vskip3pt
\centerline{$\frac{\partial^2u}{\partial xy}=0,
\ \ \text{that is}, \ \
\big(\frac{\partial}{\partial x}\big)\big(\frac{\partial}{\partial y}\big)u=0.$}
\vskip3pt
\noindent In this case, we consider $\H=L_2(\mathbb{R}^2)$
and the one-parameter unitary group $t\to U_t$ on $\H$,
where $(U_t(f))(x, y)=f(x+t,y)$, $f(x, y)\in \H$. From Stone's Theorem, there
is a self-adjoint operator $H$ defined on a certain domain in $\H$ such that
$U_t=\exp itH$. In this case, $iH$ can be identified as differentiation on
$\H$ (with domain $\D(H)$) with respect to $x$, that is $\frac{\partial}{\partial x}$.
Similarly, corresponding to the one-parameter unitary group $t\to U_t$ on $\H$,
where $(U_t(f))(x, y)=f(x,y+t)$, there is a self-adjoint operator $K$
defined on a certain domain in $\H$ and $iK$ can be identified as differentiation
$\frac{\partial}{\partial y}$ on $\H$ (with domain $\D(H)$). When $\frac{\partial}{\partial x}$
and $\frac{\partial}{\partial y}$ commute (which is guaranteed by continuous differentiability), the self-adjoint operators $H$ and $K$
commute and they are affiliated with some abelian von Neumann algebra $\R$. The (commuting) operators $A=\frac{\partial}{\partial x}$ and $B=\frac{\partial}{\partial y}$
are affiliated with $\R$.
From Theorem 3.5,
$[\N_A+\N_B]\subseteq\N_{A\,\hat\cdot\, B}$, and $[\N_A+\N_B]=\N_{A\,\hat\cdot\, B}$
when $\N_A\cap\N_B=\{0\}$.

\end{document}